\documentclass{article}

\usepackage{amsmath}
\usepackage{amssymb}
\usepackage{bm}
\usepackage{dsfont}
\usepackage{mathtools}
\usepackage[mathscr]{eucal}
\usepackage{bbm}
\usepackage{hyperref}
\usepackage{subcaption}
\usepackage{graphicx}
\usepackage{xcolor}

\DeclareMathOperator*{\argmax}{arg\,max}
\DeclareMathOperator*{\argmin}{arg\,min}

\newcommand{\pMh}{\widehat{p}_{\text{ptc}}}
\newcommand{\ph}{\widehat{p}_{\text{hist}}}

\graphicspath{ {./figs}{.} }

\usepackage{geometry}
\numberwithin{equation}{section}

\makeatletter
\def\thanks#1{\protected@xdef\@thanks{\@thanks
		\protect\footnotetext{#1}}}
\makeatother

\begin{document}

\title{The Poisson tensor completion parametric estimator}

\author{
	Daniel M. Dunlavy$^*$, 
	Richard B. Lehoucq$^\dagger$, 
	Carolyn D. Mayer$^\ddagger$, 
	and Arvind Prasadan$^\S$\\
	\emph{Sandia National Laboratories}, Albuquerque, NM and Livermore, CA
	\thanks{
		$^*$\texttt{dmdunla@sandia.gov},  
		$^\dagger$\texttt{rblehou@sandia.gov}, 
		$^\ddagger$\texttt{cdmayer@sandia.gov}, 
		$^\S$\texttt{aprasad@sandia.gov}
	}
}

\date{}

\maketitle

\begin{abstract}
We introduce the \emph{Poisson tensor completion} (PTC) estimator that exploits inter-sample relationships to compute a low-rank Poisson tensor decomposition of the frequency histogram for samples of a multivariate distribution.
Our crucial observation is that the histogram bins are an instance of a space partitioning of counts and thus can be identified with a spatial non-homogeneous Poisson process.
The Poisson tensor decomposition leads to a completion of the mean measure over all bins---including those containing few to no samples---and leads to our proposed estimator.
A Poisson tensor decomposition models the underlying distribution of the count data and guarantees non-negative estimated values obviating the need for additional constraints to ensure non-negativity. 
Furthermore, we demonstrate that our PTC estimator is a substantial improvement over standard histogram-based estimators for sub-Gaussian probability distributions because of the concentration of norm phenomenon.
\end{abstract}

\begin{figure}[b!]
	\centering
	\includegraphics[width=\textwidth]{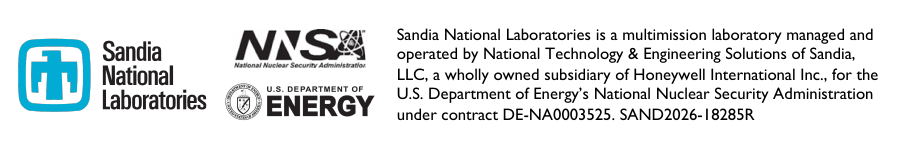}
\end{figure}

\section{Introduction}
\label{sec:introduction}

Our contribution is the \emph{Poisson tensor completion} (PTC) parametric estimator that explicitly models a multivariate distribution from samples by exploiting the relationship between frequency histograms, spatial Poisson processes, and low-rank Poisson tensor decompositions.
We demonstrate that tensor completion provides a better, at times dramatically so, representation of the histogram data. 
The PTC estimator is comprised of two steps. 
\begin{enumerate}
\item
The first step computes a low-rank Poisson tensor decomposition~\cite{chikolda12,lvdllp-arxiv:25} of the histogram bin counts to approximate the Poisson intensity measure, which is subsequently normalized to be a probability density. 
\item
The second step approximates the expectation using the  low-rank tensor approximation of the multivariate density and so is an example of a plug-in estimator.
\end{enumerate}
The first step leverages inter-sample relationships computed by a low-rank decomposition of the frequency histogram, which is an instance of a tensor, or equivalently a multi-dimensional array. 
An important consequence is that the samples conspire to enable a parametric estimator for the Poisson mean measure, which is then normalized to attain a density approximation. 
An example expectation of the second step is to approximate the differential entropy, which we use to assess the quality of the density estimator in the first step.
We focus on entropy because a Poisson tensor decomposition models the underlying distribution of the count data and guarantees a non-negative density approximation and obviates the need for additional constraints to ensure non-negativity needed for the entropy.

A recent result~\cite{han2020optimal} explains that histogram-based estimators for the differential entropy of a multivariate random variable are asymptotically optimal. 
Unfortunately, in practice, we are often limited to a finite sample of the unknown probability density upon which the differential entropy relies. 
It is also well understood that the accuracy of histogram-based probability density estimators depends upon an exponentially large number of bins in the number of variates, see, e.g.,~\cite{scott2015multivariate}. 
A kernel density estimate (KDE) is a well-known non-parametric method to estimate the underlying probability density by smoothing the data.
Imputing density values for the many empty bins that arise with an increasing number of variates is not a consequence of a KDE, which also favors nearby samples over samples further away and so is understood to be a ``local'' method.
Parametric density estimation is a global method and requires that the model is selected and a nontrivial estimation of the parameters is undertaken.
This raises the question of whether histogram data can be better leveraged by exploiting inter-sample relationships. 
Such an approach avoids the complications and leverages the strengths of parametric and KDE methods.

We believe our estimator is the first tensor-based estimator that exploits the underlying spatial non-homogeneous Poisson process related to the histogram explicitly when estimating the probability density with low-rank tensor decompositions or tensor completion.
Furthermore, we demonstrate that our PTC estimator is a substantial improvement over standard histogram-based estimators for sub-Gaussian probability distributions because of the concentration of norm phenomenon.
We also explain that heavy-tailed distributions, which do not concentrate in norm, do not benefit from the PTC estimator.
Our experiments also indicate that in practice, the number of variates is limited to available computational resources and areas for further improvements. 

Although outside the scope of this paper, the resulting tensor-completion low-rank density estimate can improve both accuracy and numerical stability in downstream high-dimensional analyses, such as hypothesis testing and point-process inference, by mitigating the severe sparsity and zero-bin pathologies that arise in large-dimensional histograms.
 
The differential entropy is a fundamental characterization  in a multivariate random variable quantifying the amount of surprise or predictability.
The entropy also has practical use in a number of applications in addition to the intrinsic mathematical and statistical desire to estimate entropy from sample observations. 
Example applications in statistical modeling include goodness-of-fit tests and tests of uniformity~\cite{gokhale1983entropy, parzen1991goodness}, alternatives to maximum likelihood estimators in settings where such estimators are not consistent~\cite{cheng1983estimating}, feature selection in biostatistics and machine learning models with tens of thousands of features~\cite{zhu2008feature}, and independent component analysis applied to blind source separation problems~\cite{hyvarinen1997new, faivishevsky2008ica}. In addition, differential entropy has been used in quantifying the thermodynamics~\cite{talukdar2017memory} of a computational process.
The entropy also has the feature that it is small where the density is small.

The remainder of our paper is organized as follows: \S\ref{sec:background} present the mathematical background and related work related to the PTC estimator, \S\ref{sec:tensor-method} provides a detailed description of the PTC estimator and an  error analysis. 
Finally, \S\ref{sec:experiments} compares and demonstrates the improvements of the PTC estimator over histogram estimators, the role of tensor rank and several practical considerations. 
Conclusions and potential future work are presented in Section~\ref{sec:conclusion}.

\section{Background: frequency histograms, spatial Poisson processes, count tensors, and related work}
\label{sec:background}

We first review the frequency histogram approximation of a $d$-variate density $p$ under which the samples are independently drawn. 
We then demonstrate that the histogram counts can be identified with a spatial Poisson process and so enable us to leverage the low-rank decomposition for count tensors introduced by Chi and Kolda in~\cite{chikolda12}.
We conclude with a discussion of related work.

\subsection{Frequency histograms}\label{ss2.1}
Suppose we have a collection of points 
\begin{align} \label{ran-pts}
x_1, \ldots, x_s \in \mathbbm{R}^d
\end{align}
independently distributed with respect to a multivariate density $p$ and then partition a finite volume 
\begin{align}\label{region-partition}
B &= \cup_{j=1}^{n_1 n_2 \cdots n_d} B_j \subset \mathbbm{R}^d \qquad  B_i \cap B_j =  \emptyset \text{ when } i \neq j
\end{align}
where $n_i$ is the number of bins along the $i$-th axis for $i=1,\ldots, d$ and denote
\begin{align*}
 n \coloneqq n_1 n_2 \cdots n_d
\end{align*}
so that $|B| = \sum_{j=1}^{n} |B_j|$ is finite.

Let $c_j$ denote the count of the points \eqref{ran-pts} that lie in bin $B_j$. 
The resulting frequency histogram can be normalized so that the approximation
\begin{align}\label{histo-q}
\ph \mathbbm{1}_{B}(x) \coloneqq  \sum_{j=1}^{n} \frac{c_j}{s |B_j|}\,\mathbbm{1}_{B_j}(x) \approx p \mathbbm{1}_{B}(x)
\end{align}
where $|B_j| = \int_{B_j} dV$ is the volume of $B_j$ and the indicator function is defined by
\begin{align*}
 \mathbbm{1}_{B_j}(x)& 
 \coloneqq \begin{cases}
 1 &  x \in B_j \\
 0 & \text{otherwise}
\end{cases}
\end{align*}
holds over $\mathbbm{R}^d$ and implies that the error depends upon the size of $\ph-p$ over the finite volume region $B $ and the size of $p$ outside of $B$.

We can then approximate the differential entropy  of random variable $\mathbf{X}\sim p$ with
\begin{subequations}
\begin{align} \label{diff-ent}
\text{ent}(p \mathbbm{1}_{B}) & = -\int_{B \mathbb{ \subset R}^d} p \log{p}\, \text{d} V 
= -\mathbb{E} \big(\log{p(\mathbf{X}) \mathbbm{1}_{B}}(\mathbf{X})\big)
\end{align}
and so
\begin{multline}\label{ent_ph-decomp}
\text{ent}(\ph\, \mathbbm{1}_{B})  = -\sum_{j=1}^{n} \int_{B_j} \ph \log \, \ph \, dV 
 = -\sum_{j=1}^{n} \frac{c_j}{s} \log{\Big(\frac{c_j}{s\,|B_j|}}\Big)
 = -\sum_{j=1}^{n} \frac{c_j}{s}\log{\frac{c_j}{s}} + \sum_{j=1}^{n} \frac{c_j}{s}\log{|B_j|} \, .
\end{multline}
Note that $\text{ent}(\ph\, \mathbbm{1}_{B})$ remains well-defined for finite bin volumes $|B_j|$ and zero counts because $\lim_{x \to 0^+} x\log{x}=0$.
\end{subequations}

\subsection{Spatial Poisson processes}\label{ss2.2}

We now show that a model for the random set of points \eqref{ran-pts} is a spatial Poisson process.
Define the counting measure 
\begin{align*}
N(A) &  \coloneqq \sum_{k=1}^{s}\mathbbm{1}_{A}(x_k)
\end{align*}
that returns the number of points \eqref{ran-pts} in  $A \subset \mathbb{R}^d$
so that 
\begin{align*}
N(B_j) = \sum_{k=1}^{s}\mathbbm{1}_{B_j}(x_k) = c_j\,.
\end{align*}
In words, the number of points in the frequency histogram bin $B_j$ is given by $N(B_j) $.
Well-known (see, e.g., \cite[p.53]{scott2015multivariate}) is that the count 
\begin{subequations}\label{P-approx}
\begin{align} \label{suc-prob}
c_j \sim \text{Binom}\big(s,  \mu(B_j) \big) \text{ where } \mu(B_j) = \int_{B_j} p \, \text{d}V\,,
\end{align}
i.e., is binomially distributed with $s$ events and success probability $\mu(B_j)$.
Approximation of the Binomial distribution by the Poisson distribution is excellent when $s$ is large and $\mu(B_j)$ is small so that 
\begin{align}\label{mm}
\nu_j = s\,\mu(B_j)
\end{align}
is moderate.
\end{subequations}
The bin counts $c_j$ can then be approximated with a spatial (non-homogenous) Poisson process because the following two properties are satisfied:
\begin{enumerate}
\begin{subequations} \label{spp-cond}
\item
For each bin $B_j$,
\begin{align}
\mathbb{P} \lbrack N(B_j) = k \rbrack = 
e^{-\nu_j}\, \frac{\nu_j^k}{k\, ! } 
\end{align}
where the success probability \eqref{suc-prob} is proportional to the mean measure \eqref{mm} for the spatial Poisson process and $p$ can be identified with the Poisson intensity. 
This identification is the crucial relationship between the multivariate random variable density $p$, the Poisson intensity measure, and the spatial Poisson process that are leveraged by the PTC estimator.
\item
Given the disjoint bins $B_1\,, \ldots, B_{n}$, then 
\begin{align}
N(B_1)\,, \,\cdots, N(B_{n}) \label{spp-condb}
\end{align}
are independent random variables.
This holds  because we assumed that the set of points \eqref{ran-pts} are randomly drawn under $p$ and we remark that the independence among the $N(B_j)$ holds regardless of the independence of variates.
\end{subequations}
\end{enumerate}
We replace the histogram estimate $\ph$ of the density $p$ with the approximation $\pMh$ to be described \S\ref{sec:tensor-method} and then estimate the entropy as in \eqref{ent_ph-decomp}.
Our estimator is parametric in the $n = n_1 \cdots n_d $ mean measures $\nu_j$. 
A subsequent normalization renders an approximation of the density $p$ over each bin, which under mild regularity assumptions converges when the maximum bin size decreases to zero.

\subsection{Count tensors} \label{ss2.3}
A tensor is multi-dimensional array where the order is the number of dimensions $d$, which for our application is the number of variates.
Examples include scalars (tensors of order zero), e.g., $x$; vectors (tensors of order one), e.g., $\mathbf{t}$; and matrices (tensors of order two), e.g., $\mathbf{T}$. 
Following tensor notation in~\cite{KoBa09}, we denote tensors of order three and higher with bold capital script letters, e.g., ${\boldsymbol{\mathscr{T}}}$.

We define the \emph{histogram tensor} for a $d$-variate density $p$ containing the counts $c_j$ defined in \S\ref{ss2.1} using the partitioning scheme specified in \eqref{region-partition} as 
\begin{align*}\label{hist-ten}
\boldsymbol{\mathscr{T}} \in \mathbb{Z}_{+}^{n_1 \times n_2 \times \cdots \times n_d}\,,
\end{align*}
where $n_i$ is the number of bins into which dimension $i$ is partitioned and $\mathbb{Z}_{+}$ is the set of non-negative integers. 
Unless the number of points $s$ is exponentially large in the number of variates, the histogram-based density estimator $\ph$ is sparse with many zero or near-zero bins, and so does not provide an accurate approximation of the probability density needed to estimate an expectation, e.g., entropy.
We exploit histogram sparsity by computing a low-rank Poisson canonical polyadic (CP) tensor decomposition in~\cite{chikolda12} of ${\boldsymbol{\mathscr{T}}}$ and use tensor completion to provide a model of the expected counts in all bins including those containing no sample points. 
A key assumption made by Chi and Kolda in~\cite{chikolda12} is the independence of  $n=n_1 \cdots n_d$ Poisson random variables, which is satisfied by the $n$ counting measures $N(B_1), \cdots, N(B_n)$. 
These counting measures determine the mean measures $\nu_1, \cdots, \nu_n$ approximated via the solution of a maximum likelihood problem. Although Chi and Kolda introduced the low-rank Poisson CP tensor decomposition in~\cite{chikolda12}, they did not make a connection between their work and multivariate distributions. Our work here is the first to leverage the connections between frequency histograms of samples for multivariate distributions, a spatial non-homogeneous Poisson process, and their work on Poisson tensor decompositions.

The selection role of tensor rank


We develop our estimator using the Poisson CP tensor decomposition, however 
other low-rank tensor Poisson decompositions---e.g., Tucker~\cite{tucker} or tensor train~\cite{tensor_train} decompositions---when such models and associated statistical theory are adapted for use with Poisson data.

\subsection{Related work}
\label{sec:related}

The authors of \cite{vand:21} introduce the use of tensor decompositions as an improvement on the histogram and KDE for high-dimensional density estimation. 
Three other papers \cite{amks:21, amks:22a, amks:22b} also use tensor decompositions to estimate a multivariate density. 
The report \cite{amks:22a} is similar in spirit to that in \cite{vand:21}, where a transformation into the Fourier domain leads to low-rank structure assumed in \cite{vand:20, vand:21}.  
The recent report \cite{xppw:25} estimates multivariate intensity functions of spatial point processes using matrix- and tensor-based methods. 
We note that our estimator is distinguished from these tensor-based estimators in that we model the histogram tensor values as counts and compute low-rank Poisson tensor decompositions to satisfy the histogram modeling assumptions introduced in \S\ref{ss2.2} and \S\ref{ss2.3}. Specifically, the Poisson CP decomposition of the histogram tensor models the mean Poisson intensity measures across all bins associated with the spatial Poisson process described in \S\ref{ss2.2}; this is discussed in detail in \S\ref{sec:tensor-method}. Moreover, the approaches in \cite{vand:21,xppw:25}, which are the most closely related to our PTC estimator, rely on tensor decompositions fit by minimizing squared error loss, which may lead to violations of model assumptions and poor model fits for count data as discussed in \cite{chikolda12}.

The tensor-based estimators in \cite{vand:20, vand:21}, i.e., multi-view models, have also been studied in the context of latent variable models, where parameter inference is performed using the method of moments \cite{Anandkumar14}. As we focus here on parameter inference uising maximum likelihood estimation as noted in \S\ref{ss2.3} direct comparisons with that work is beyond the scope of this paper and left as future work.

We reviewed in \S\ref{ss2.1} that the differential entropy of a random variable is the expectation of the negative logarithm of corresponding probability density. 
Classical approaches to estimating the multivariate probability density include histograms and kernel density estimators (KDE) with the limitation that given $s$ samples in $d$ dimensions, consistency generally requires $s h^d \rightarrow \infty$, where $h$ is the bin-width of the histogram or the bandwidth parameter for the KDE, see, e.g., the textbook \cite{scott2015multivariate} for a review.
For a general survey of non-parametric entropy estimation, we refer the interested reader to \cite{beirlant1997nonparametric}. 

The accuracy of the histogram-based differential entropy estimator \eqref{ent_ph-decomp} depends upon an exponentially large number of bins with an increasing number of variates.
The plug-in estimator 
\begin{equation}\label{eq:plugin_estimator}
	-\frac{1}{s} \sum_{i} \log \ph(x_i)\,
\end{equation}
is simpler, however, its accuracy is also beholden to an exponential number of samples \cite{joe1989estimation}.  Our work directly improves upon the plug-in estimator \eqref{eq:plugin_estimator} that is obtained with a histogram or KDE density estimate for a fixed sample size.
A slight improvement to \eqref{eq:plugin_estimator} might be obtained by splitting the data, that is, by using $s - 1$ points to estimate the density and then using the held-out point to estimate the entropy, before averaging over all held-out points, as in \cite{hall1993estimation}.

In one dimension, the inverse of the spacing between samples yields a rough estimate of the density at a given point, which can be used (with some bias-correction terms) in \eqref{eq:plugin_estimator} using the scheme in \cite{tarasenko1968evaluation}. 
In larger dimensions, the natural analog of the spacing estimator is based on the $k$-nearest neighbor distances. 
The work in \cite{kozachenko1987sample} is such an estimator, which uses $1$-nearest neighbor distances; later work, including \cite{singh2003nearest} and \cite{goria2005new} generalized the estimator to $k$-nearest neighbor distances for $k > 1$.
We refer the interested reader to \cite{berrett2019efficient, hall1984limit,  mnatsakanov2008k} for a more thorough characterization of the spacing and $k$-NN estimators, respectively. 
In short, for $s$ samples and $k$ a function of $s$, we require $k$ to be chosen so that $k(s) / s \rightarrow 0$ for the estimator to be asymptotically unbiased \cite[Lemma~3]{berrett2019efficient}.
At the same time, \cite[Lemma~7]{berrett2019efficient} indicates that $k(s) \rightarrow \infty$ is required for a function of $s$ for the variance to be minimized, i.e., for the estimator to be efficient. 
In general, a choice of $k_s \propto \log^6 s$ is sufficient for both conditions to hold \cite[Theorem~1]{berrett2019efficient}. 
As a point of interest, the requirement that a density have at least one continuous derivative at all points (see  \cite[Theorem~1]{berrett2019efficient}) precludes the application of the above theory to the uniform distribution. A recent density estimator of non-trivial complexity is that in \cite{ariel2020estimating}, which seeks to estimate a set of marginal distributions and an accompanying  copula or dependence structure of said marginals. 
The performance of this estimator is generally worse than the $k$-nearest neighbor ($k$-NN) distances \cite{mack1979multivariate}, which is easily adapted to yield an estimator for the entropy and is based upon pairwise distances between samples. As we show in Section~\ref{sec:experiments}, our new PTC differential entropy estimator improves upon estimators based on these $k$-NN density approximations for several classes of probability distributions.

\section{Poisson tensor completion}
\label{sec:tensor-method}

We reviewed in \S\ref{ss2.1}--\S\ref{ss2.3} the relationship between a frequency histogram, a spatial Poisson process modeling the sample points drawn from the $d$-variate distribution defined by density $p$, and an order-$d$ tensor of counts $ \boldsymbol{\mathscr{T}}$.
In this section, we leverage these relationships to define the \emph{Poisson tensor completion} (PTC) estimator and demonstrate how it can be used for efficient estimation. 

\subsection{Density estimation}

Tensor completion is a method for computing a low-rank model of tensor data given a sample of that data; see \cite{tensor_completion} and the references therein for an introduction to tensor completion and a survey of the many approaches developed. A standard approach to tensor completion is to compute a low-rank decomposition of the sampled data and then estimate the unobserved values in the tensor from that low-rank model. Here, we use that approach and compute a Poisson canonical polyadic (CP) tensor decomposition of the histogram tensor, an approach novel with this paper.

Recall that in \S\ref{ss2.3} we defined the histogram tensor $\boldsymbol{\mathscr{T}}$ for a $d$-variate density $p$.
Using multi-index notation to indicate a tensor element, i.e., $t_\mathbf{i} \equiv t_{i_1 i_2 \cdots i_d}$ denotes entry $ \mathbf{i}=(i_1,\dots,i_d)\in [n_1] \otimes [n_2] \otimes \dots \otimes [n_d]$ of $\boldsymbol{\mathscr{T}}$ with $[n_i] = \{1, \dots, n_i\}$, we follow Chi and Kolda~\cite{chikolda12} and model the elements of the histogram tensor as independent Poisson random variables so that
\begin{align}\label{eq:poisson}
	t_{\mathbf{i}} \sim \mbox{Poisson}\left(m_{\mathbf{i}}\right)\,.
\end{align}
In other words, $t_\mathbf{i} $ is the number of samples $x_i$ that are located in bin $B_\ell$ for a unique linear index $\ell \in [n]$ so that the Poisson distribution assumption~\eqref{eq:poisson} is equivalent to the spatial Poisson process~\eqref{spp-cond} modeling the bin counts. 
The linear index $\ell$ denotes a mapping between bin and tensor multi-indices and  an example is the natural ordering index notation for tensors introduced in~\cite{Ballard_Kolda_2025} where the relationship is $\ell = \mathbb{L}(\mathbf{i})$ and $\mathbf{i} = \mathbb{T}(\ell)$ and $\mathbb{L}:[n_1] \otimes \cdots \otimes [n_d] \mapsto \mathbb{Z}_+$ and $\mathbb{T}:\mathbb{Z}_+ \mapsto [n_1] \otimes \cdots \otimes [n_d]$.

To the best of our knowledge, we are the first to identify the relationship between a spatial Poisson process and a Poisson canonical polyadic (CP) tensor decomposition~\cite{chikolda12} to model histogram bin counts. 
We also remark that $\boldsymbol{\mathscr{T}}$ will contain many entries corresponding to unobserved or low counts unless the number of samples $s$ is exponentially large in the number of dimensions (i.e. variates) $d$. 
Poisson tensor completion is a mechanism for imputing the expected values for these entries under the Poisson assumption~\eqref{eq:poisson}, i.e., the spatial Poisson process~\eqref{spp-cond} modeling the bin counts.

We  complete $\boldsymbol{\mathscr{T}}$ by imposing low-rank CP structure on the corresponding Poisson parameter tensor
\begin{align} \label{M-def}
	 \boldsymbol{\mathscr{M}} \coloneqq \sum_{r=1}^R \lambda_r \, \mathbf{a}_r^{(1)} \circ \mathbf{a}_r^{(2)} \circ \cdots \circ \mathbf{a}_r^{(d)} \in \mathbb{R}_{+}^{n_1 \times n_2 \times \cdots \times n_d}\; ,
\end{align}
where $\circ$ denotes the outer product of such that $\|\mathbf{a}_r^{(i)}\|_1 = 1 \, , \forall i \in \{1,\dots,d\}, \forall r \in \{1,\dots,R\}$, 
and $\mathbb{R}_{+}$ is the set of positive real values. 
Note that the entries $m_\mathbf{i}$ of the low-rank CP parameter tensor $\boldsymbol{\mathscr{M}}$ are positive by definition and they model the probability of a number of events occurring in a bounded region of $\mathbb{R}^d$ as defined by the spatial Poisson process associated with the frequency histogram discussed in \S\ref{sec:background}.
 
Let $ m_{\mathbf{i}}$ is the $\mathbf{i}$ entry of the rank $R$ tensor $\boldsymbol{\mathscr{M}}$. We estimate the mean measure~\eqref{mm} of the spatial Poisson process using the maximum Poisson likelihood estimation 
 \begin{equation}\label{ml-est}
 \begin{aligned}\boldsymbol{\widehat{\mathscr{M}}} &= \argmax \sum_{\mathbf{i} = (1, 1, \cdots,1)}^{(n_1,n_2,  \cdots, n_d)}( t_\mathbf{i}\log m_{\mathbf{i}} - m_{\mathbf{i}}) &= \argmin \sum_{\mathbf{i} = (1, 1, \cdots,1)}^{(n_1,n_2,  \cdots, n_d)}( m_{\mathbf{i}} -t_\mathbf{i}\log m_{\mathbf{i}} )
 \end{aligned}
 \end{equation}
 as introduced in~\cite{chikolda12}.
Then, with the natural ordering index mapping $\mathbb{T}$ defined following \eqref{eq:poisson}, each entry in $\boldsymbol{\widehat{\mathscr{M}}}$ estimates the true Poisson parameter of~\eqref{eq:poisson}, which in turn represents the expected count of the corresponding bin in the frequency histogram as follows:
 \begin{align} \label{mhat-lilkeli}
 	 \widehat{m}_{\mathbb{T}(\ell)} = \widehat{m}_{\mathbf{i}}  \approx m_\mathbf{i} = \mathbb{E}\left(t_\mathbf{i}\right) \, .
\end{align}
We can therefore denote the induced probability density
\begin{align}\label{pMh}
\pMh \, \mathbbm{1}_B(x) & \coloneqq \frac{1}{\| \boldsymbol{\widehat{\mathscr{M}}}\|_1} \sum_{j=1}^{n} \frac{\widehat{m}_{\mathbb{T}(j)}}{|B_j|} \, \mathbbm{1}_{B_j} (x)
\end{align}
where $\|\boldsymbol{\widehat{\mathscr{M}}}\|_1 = \sum_{r=1}^R \widehat{\lambda}_r $ because the vectors of a Poisson CP tensor decomposition \eqref{M-def} are normalized to unit length in the one-norm. 
We remark that when the rank $R$ is greater than one, the rank one tensors comprising the summands in~\eqref{M-def} model dependencies among the variates, and unit rank implies that the variates are independent. 
A striking example is a mixture model where the rank models the number of components in a mixture model so that $\lambda_r/\sum_{r=1}^R \widehat{\lambda}_r$ is the weight of the $r$-component; see \S\ref{s:gms} for a discussion and experiments.

Our approximation is parametric because 
\begin{align*}
\widehat{m}_{\mathbf{i}} \approx  \nu_{\mathbf{i}} = s\int_{B_{\mathbf{i}}} p \, \text{d}V\,,
\end{align*}
i.e., each entry approximates the mean measure over a bin.
Our estimate of the underlying multivariate density is given by a normalization of the dense tensor $\boldsymbol{\widehat{\mathscr{M}}}$  that is the same size as the histogram tensor $\boldsymbol{\mathscr{T}}$, which contains $\prod_{i=1}^d n_i$ entries and subject to the available computation resources, may be too large to store. 
Thus, we store $\boldsymbol{\widehat{\mathscr{M}}}$ in its decomposed form defined by~\eqref{M-def}, which is comprised of $(R+1)\sum_{i=1}^d n_i$ entries, computing $\|\boldsymbol{\widehat{\mathscr{M}}}\|_1$ and $\widehat{m}_{\mathbb{T}(j)}$ in~\eqref{ptc-est} only as needed. 
Computing each $\widehat{m}_{\mathbb{T}(j)}$ requires $R(d+1)$ operations.

Our introduction raised the question of whether histogram data can be better leveraged by exploiting inter-sample relationships inherent in the Poisson likelihood problem \eqref{ml-est}, which we now discuss.
Inserting \eqref{M-def} into the $\mathbf{i} \equiv (i_1, i_2, \cdots, i_d)$-th entry of the sum in \eqref{ml-est} grants
\begin{multline*}
m_{\mathbf{i}} -t_\mathbf{i}\log m_{\mathbf{i}} = \sum_{r=1}^R \lambda_r \, \mathbf{a}_r^{(i_1)} \circ \mathbf{a}_r^{(i_2)} \circ \cdots \circ \mathbf{a}_r^{(i_d)} 
- t_{(i_1,\dots,i_d)}\log \Big(\sum_{r=1}^R \lambda_r \, \mathbf{a}_r^{(i_1)} \circ \mathbf{a}_r^{(i_2)} \circ \cdots \circ \mathbf{a}_r^{(i_d)} \Big)
\end{multline*}
where $\sum_{r=1}^R \lambda_r\mathbf{a}_r^{(i_1)} \circ \mathbf{a}_r^{(i_2)} \circ \cdots \circ \mathbf{a}_r^{(i_d)}$ is the $(i_1, i_2, \cdots, i_d)$ entry of $ \boldsymbol{\mathscr{M}} $.
This explains that each count $t_{(i_1,\dots,i_d)}$ is directly involved in the determination of the $ n_k $ elements in the mode-$k$ fiber\footnote{The mode-$k$ fiber of a tensor is the vector that results when all indices but the $k$-th are held fixed.} for $k=1, \cdots, d$
of the rank $R$ tensor $ \boldsymbol{\mathscr{M}} $. 
Because there are $n_1 n_2 \cdots n_d$ counts, then the ratio of the total number of counts to the number of parameters satisfies
\begin{align*}
\frac{n_1 n_2 \cdots n_d}{R(n_1 + n_2 + \cdots + n_d)}  \to \infty
\end{align*}
as the number of variates increases.
This implies that an exponential amount of count data is used to determine a linear number of elements in $ \boldsymbol{\mathscr{M}} $ and explains how each $t_{(i_1,\dots,i_d)}$ effects elements of $ \boldsymbol{\mathscr{M}} $.

The recent report \cite{lvdllp-arxiv:25} demonstrates that the likelihood problem \eqref{ml-est} is well-posed and explains the role of the rank $R$ in identifiability and indeterminacy.
The authors of \cite[\S2.1]{lvdllp-arxiv:25} also provide a necessary and sufficient condition on the count tensor $\boldsymbol{\mathscr{T}}$ so that all the elements of $\boldsymbol{\widehat{\mathscr{M}}}$ are positive. In practice, we leverage the Poisson CP decomposition approach of \cite{chikolda12}, which guarantees $\boldsymbol{\widehat{\mathscr{M}}}$ to be non-negative for well-posed inference problems. 
We demonstrate in \S\ref{tensor-thres} that the appearance of small positive numbers and zeros in $\boldsymbol{\widehat{\mathscr{M}}}$ can be beneficial in computing PTC estimators efficiently.

\subsection{Plug-in entropy estimation}

The expected value of $f(X)$ where $X\sim p$ for a real-valued function $f: \mathbb{R}^d \to \mathbb{R}$ is
\begin{align}\label{exp-value}
 \int_{\mathbb{R}^d} f \pMh\, \mathbbm{1}_B \, dV & = \int_B f  \pMh\, dV
                               = \sum_{j=1}^{n} \int_{B_j} f\, \pMh\, dV 
                               = \frac{1}{\|\boldsymbol{\widehat{\mathscr{M}}}\|_1} \sum_{j=1}^{n} \widehat{m}_{\mathbb{T}(j)} \, \bar{f}_{B_j} 
\end{align}
where $\bar{f}_{B_j} \coloneqq \frac{1}{|B_j|}\int_{B_j} f\, dV$ is the average value of $f$ over the bin $B_j$.

Setting $f = \log (\pMh\, \mathbbm{1}_B)$ defines the \emph{PTC estimator for differential entropy} as 
\begin{align}
\label{ptc-est}
\text{ent}(\pMh\, \mathbbm{1}_B) & = -\sum_{j=1}^{n} \frac{\widehat{m}_{\mathbb{T}(j)}} {\|\boldsymbol{\widehat{\mathscr{M}}}\|_1}\log\left(\frac{\widehat{m}_{\mathbb{T}(j)}}{\|\boldsymbol{\widehat{\mathscr{M}}}\|_1 |B_j|}\right)
\end{align}
because $\bar{f}_{B_j} = \log\left(\frac{\widehat{m}_{\mathbb{T}(j)}}{\|\boldsymbol{\widehat{\mathscr{M}}}\|_1 |B_j|}\right)$.

\subsection{Error Analysis}\label{subsec:error-analysis}

A detailed error analysis is outside the scope of this paper. 
Instead we present a preliminary error analysis using recent analysis compelling the PTC approximation of the density $p$.

We explained at the end of \S\ref{ss2.2}  that when the Poisson approximation \eqref{P-approx}  holds, we have that 
\begin{subequations}
\begin{align}\label{phat}
\widehat{p}\,\mathbbm{1}_{B} \coloneqq \sum_{j=1}^n \frac{\nu_j}{s|B_j|}  \mathbbm{1}_{B_j} \to p \mathbbm{1}_{B} \text{ as } n \to \infty
\end{align}
where we assumed that $|B| = \sum_{j=1}^{n} |B_j|$ is finite in \S\ref{ss2.1} and a calculation using \eqref{mm}  establishes  that
\begin{align}
\int_B \widehat{p}\,\text{d}V= \int_B p \,\text{d}V\,.
\end{align}
\end{subequations}
Let $\boldsymbol{\mathscr{N}}$ denote the tensor containing the $n$ mean measures $\nu_j$ so that $\boldsymbol{\mathscr{N}}$ is ordered consistently with $\boldsymbol{\widehat{\mathscr{M}}}$.

The recent report \cite[Thm~4]{lplvld:25} builds upon the results in \cite{lvdllp-arxiv:25} to show that
\begin{align}
\| \boldsymbol{\widetilde{\mathscr{M}}} - \boldsymbol{\mathscr{N}}\|_F^2 
\leqslant C \frac{\sqrt[d]{n}\,(d+\delta)^2 d \,R^2}{(\beta + 1)^2}
\end{align}
holds for $\delta > 0$  with probability no larger than $1-\delta$
for a constant $C$ where we assume that $k = n_1 = \cdots = n_d$ so that $k=\sqrt[d]{n}$, $\beta$ is a positive lower bound over all $n=k^d$ Poisson parameters over the partition \eqref{region-partition} of $B$, and
\begin{align}
\label{llMt}
\boldsymbol{\widetilde{\mathscr{M}}} = \argmax \sum_{\mathbf{i} = (1, 1, \cdots,1)}^{(k,k,  \cdots, k)}\big( (t_\mathbf{i}+1)\log (m_{\mathbf{i}}+1) - (m_{\mathbf{i}}+1) \big)\,.
\end{align}
The authors of \cite{lplvld:25} used $\boldsymbol{\widetilde{\mathscr{M}}}$ to avoid zero counts in $\boldsymbol{{\mathscr{T}}}$ and simplify the analysis in providing an upper bound on the error.
The error in the approximate Poisson parameters $ \boldsymbol{\widehat{\mathscr{M}}} $ over the region $B$ is then
\begin{align*}
\| \boldsymbol{\widehat{\mathscr{M}}} - \boldsymbol{\mathscr{N}}\|_F \leqslant \| \boldsymbol{\widehat{\mathscr{M}}} - \boldsymbol{\widetilde{\mathscr{M}}}\|_F + \| \boldsymbol{\widetilde{\mathscr{M}}} - \boldsymbol{\mathscr{N}}\|_F 
\leqslant \| \boldsymbol{\widetilde{\mathscr{M}}} - \boldsymbol{\mathscr{N}}\|_F\,.
\end{align*}
The last inequality follows because the gradient of the likelihoods \eqref{ml-est} and \eqref{llMt} are the same.
This implies that the error over $\mathbb{R}^d$ is the above sum and the error on $\mathbb{R}^d\setminus B$.
Our assumption on $\beta$ implies that
\begin{align} \label{lbN}
\beta n \leqslant \| \boldsymbol{\mathscr{N}} \|_1 \leqslant \sqrt{n} \| \boldsymbol{\mathscr{N}} \|_F
\end{align}
so that we can conclude that the relative error in the PTC estimate of the density $p$ is
\begin{align} \label{rel-errN}
\frac{\| \boldsymbol{\widehat{\mathscr{M}}} - \boldsymbol{\mathscr{N}}\|_F}{\| \boldsymbol{\mathscr{N}} \|_F}
\leqslant  C \frac{1}{n^{1/2-1/d}} \frac{(d+\delta)^2 d \,R^2}{\beta(\beta + 1)^2} \,.
\end{align}
The relative error decreases to zero with increasing $n$ and $d>2$ and increases with $R^2$ and is inversely proportional to $\beta$.
This indicates that the PTC mean measures $\widehat{m}_{\mathbb{T}(j)}$ are converging to the mean measures $\nu_j$ so that we may formally conclude that 
\begin{align*}
\| (\pMh -  \widehat{p})\mathbbm{1}_{B} \| \to 0 \text{ as } n \to \infty \text{ for } d> 2
\end{align*} 
where $\| \cdot \|$ is an appropriate norm and $s$ is sufficiently large so that the Poisson approximation \eqref{P-approx} holds.
We may now conclude that for $d>2$ variates, 
\begin{multline}
\| \pMh -  p \|  \leqslant \| (\pMh -  \widehat{p})\mathbbm{1}_{B} \| + \|(\widehat{p} -  p)\mathbbm{1}_{B} \| 
 + \| p \mathbbm{1}_{\mathbb{R}^d\setminus B} \| \to \| p \mathbbm{1}_{\mathbb{R}^d\setminus B} \| 
\text{ as } n \to \infty
\end{multline}
and $s$ is sufficiently large so that the Poisson approximation \eqref{P-approx} holds. 

Several remarks are in order. The inequality \eqref{lbN} is instrumental for relative error \eqref{rel-errN} and use of $\ph$ does not lead to a similar bound as in \eqref {lbN} because of the exponential growth of zero bins. 
This is a compelling advantage of PTC over a histogram when $\| p \mathbbm{1}_{\mathbb{R}^d\setminus B} \|$ is small.
The relative error \eqref{rel-errN} increases quadratically with the rank $R$ indicating that the rank needed cannot be large relative to the number of variates $d$. 
The report \cite{lvdllp-arxiv:25} suggests that when the $d$ variates are moderately to strongly dependent, then $R$ is not small. 
Our final remark is that the relative error \eqref{rel-errN} only considers the error on the finite volume $B$ and so the quantity $\| p \mathbbm{1}_{\mathbb{R}^d\setminus B} \|$ is ignored.

For what class of distributions $p$ can we expect $\| p \mathbbm{1}_{\mathbb{R}^d\setminus B} \|$ to be small?
A sufficient condition is that a majority of the probability of a multivariate distribution is contained by the finite volume $B$.
In other words, the tails decay sufficiently fast.
Such a trend, as it turns out, is generic for the important class  of sub-Gaussian distributions as we now review.
The random vector $\mathbf{X} \in \mathbb{R}^d $ has sub-Gaussian components $X_i$ when
\begin{subequations}\label{sg-dists}
\begin{align}
\mathbb{P} \{ |X_i | \geqslant t \} \leqslant 2 e^{-c t^2}
\end{align}
for a constant $c$ and $i=1,\ldots, d$. In words, the tails of the multivariate distribution for $X$ decay as for a multivariate normal distribution. 
Examples of sub-Gaussian distributions include Gaussian, uniform, bounded distributions and mixtures of sub-Gaussian distributions.
An important consequence is that for such distributions there is a concentration of norm, i.e., 
\begin{align}
\mathbb{P} \{ \| \mathbf{X} \|_2   - \sqrt{d} \geqslant t \}  \leqslant 2 e^{-\tilde{c} t^2}
\end{align}
\end{subequations}
where $ \| \mathbf{X} \|_2 = \sqrt{\sum_{i=1}^d X_i^2}$ and a constant $\tilde{c}$. 
In particular, a sample of a multivariate standard normal concentrates about a thin spherical shell of radius $\sqrt{d}$ so that as the 
number of variates $d$ increases so does the distance of the spherical shell from the origin. 
This suggests a threshold on the size of the elements of $\boldsymbol{\widehat{\mathscr{M}}}$  during the computation of $\boldsymbol{\widehat{\mathscr{M}}}$ and such a scheme is introduced \S\ref{tensor-thres}.

The error in the expected value \eqref{exp-value} of $f(X)$ where $X\sim p$  is
\begin{align}
 \int_{\mathbb{R}^d} f p \,  \text{d}V -  \int_{\mathbb{R}^d} f \, \pMh\, \mathbbm{1}_B \, \text{d}V = \int_{\mathbb{R}^d} f \, (p - \pMh) \, \mathbbm{1}_B \, \text{d}V + \int_{\mathbb{R}^d} f p \mathbbm{1}_{\mathbb{R}^d\setminus B}\,  \text{d}V\,.
\end{align}
The error depends upon the approximation of $p$ by $\pMh$ in the direction of $f$ over the region $B$ and the product $f p$ over the region outside of $B$.
Our PTC convergence analysis provides conditions under which $\pMh$ converges, which implies sufficient conditions for the consistency of the PTC plug-in estimator. 
For instance, asymptotic consistency demands that as the number of samples $s$ increases that
\begin{align*}
\int_{\mathbb{R}^d} f \, (p - \pMh) \, \mathbbm{1}_B \, \text{d}V \to 0 
\text{ and } \int_{\mathbb{R}^d} f p \mathbbm{1}_{\mathbb{R}^d\setminus B}\,  \text{d}V \to 0
\end{align*}
both in probability.
A complicating factor is that the volume $B$ depends upon the samples and so a sufficient property is that the volume $B$ grows slowly with $s$.
Sub-Gaussian distributions satisfy such a property.

The error of the differential entropy of a sub-Gaussian distribution is a motivating expectation because the entropy is small where the density is small:
\begin{subequations}
\begin{multline}
\text{ent}(\pMh) - \text{ent}(p)  
=
-\sum_{j=1}^n \Bigg(  \frac{\widehat{m}_{\mathbb{T}(j)}} {\|\boldsymbol{\widehat{\mathscr{M}}}\|_1}\log\Big(\frac{\widehat{m}_{\mathbb{T}(j)}}{\|\boldsymbol{\widehat{\mathscr{M}}}\|_1 |B_j|}\Big)
- \int_{B_j} p\mathbbm{1}_{B_j} \log(p \mathbbm{1}_{B_j}) \, \text{d} V \Bigg)
 -   \text{ent}(p\, \mathbbm{1}_{\mathbb{R}^d\setminus B})
\end{multline}
and the analogous error using the histogram estimator is
\begin{multline}
\text{ent}(\ph) - \text{ent}(p)    = 
-\sum_{j=1}^n \Bigg(  \frac{c_j}{s }\,\mathbbm{1}_{B_j} \log\Big( \frac{c_j}{s |B_j|}\,\mathbbm{1}_{B_j}\Big)
- \int_{B_j} p\mathbbm{1}_{B_j} \log(p \mathbbm{1}_{B_j}) \, \text{d} V \Bigg)
-   \text{ent}(p\, \mathbbm{1}_{\mathbb{R}^d\setminus B}) \,.
\end{multline}
\end{subequations}
Both estimates enables us to conclude that error can be no smaller than the size of $p$ on $\mathbb{R}^d\setminus B$ and are consistent estimators when the entropy on the tails of the distribution converges to zero in probability. 
The estimators differ in their bin-wise expressions and consists of the error in density estimation and then entropy estimation.
The value of $\text{ent}(\ph)$ over bins without counts is zero whereas $\text{ent}(\pMh)$ is positive because of tensor completion.

\section{Experiments}
\label{sec:experiments}

In our experiments, we use differential entropy as a way to assess the PTC estimator of the density for several multivariate distributions. In  \S\ref{s:bs}, we investigate the effect of bin size on the resulting entropy approximation using the PTC estimator or a standard histogram-based estimate. We also compare the estimated entropy using PTC, histograms, and the $k$-NN method with different numbers of samples. In \S\ref{s:gms} we use a Gaussian mixture model to understand the role of tensor rank selection for the PTC estimator. 
Because the entropy is small where the density is small, we introduce a simple thresholding in \S\ref{tensor-thres} on the size of the entries of $\boldsymbol{\widehat{\mathscr{M}}}$ to reduce the cost,
which becomes prohibitive with an increasing number of variates $d$  because $\boldsymbol{\widehat{\mathscr{M}}}$ cannot be formed without a significant amount of memory nor can we compute the individual elements $\widehat{m}_{\mathbb{T}(j)}$ without a significant amount of computation. In  \S\ref{s:real-data}, we apply our PTC estimator to approximate the differential entropy for a couple real datasets.

Our experiments were performed in Python and are summarized by the following three steps.
\begin{enumerate}
\item
Sample the appropriate distribution to generate $s$ random points \eqref{ran-pts}.
\item 
Bin the random points according to a specified binning scheme. Use a sparse representation of the histogram, keeping track of the nonempty bins, counts, and bin volumes. 
Unless otherwise specified:
\begin{itemize}
\item 
for estimates of density made using a histogram directly, the histogram bins are width 
\begin{align} \label{opt-bin-size}
3.5 s^{-\frac{1}{d + 2}} \text{ in each dimension so that } |B_j| = (3.5)^d s^{-\frac{d}{d + 2}}
\end{align} 
where $s$ is the number of samples from a $d$-dimensional distribution.
We selected this value to be close to the asymptotically optimal bin widths for the multivariate Gaussian $N(0_d,I_d)$, \cite[Eq.~(3.66)]{scott2015multivariate}). 
The bin edges can be found using \texttt{arange} function within \texttt{numpy} \cite{numpy}.
\item 
histogram tensors have $n_i=20$ in each dimension for a total of $n=20^d$ bins. 
The bin edges can be found using the \texttt{histogram\_bin\_edges} function within \texttt{numpy}.
\end{itemize}
\item 
Use the routines \texttt{sptensor.from\_data} and \texttt{cp\_apr}  in the Tensor Toolbox \texttt{pyttb} \cite{pyttb} to compute the low-rank Poisson CP tensor model $\boldsymbol{\widehat{\mathscr{M}}}$ as in \eqref{M-def}. 
Unless otherwise specified, the experiments used rank $R$ at most $5$. 
 \end{enumerate}

\begin{figure}
\centering
\includegraphics[width=\linewidth]{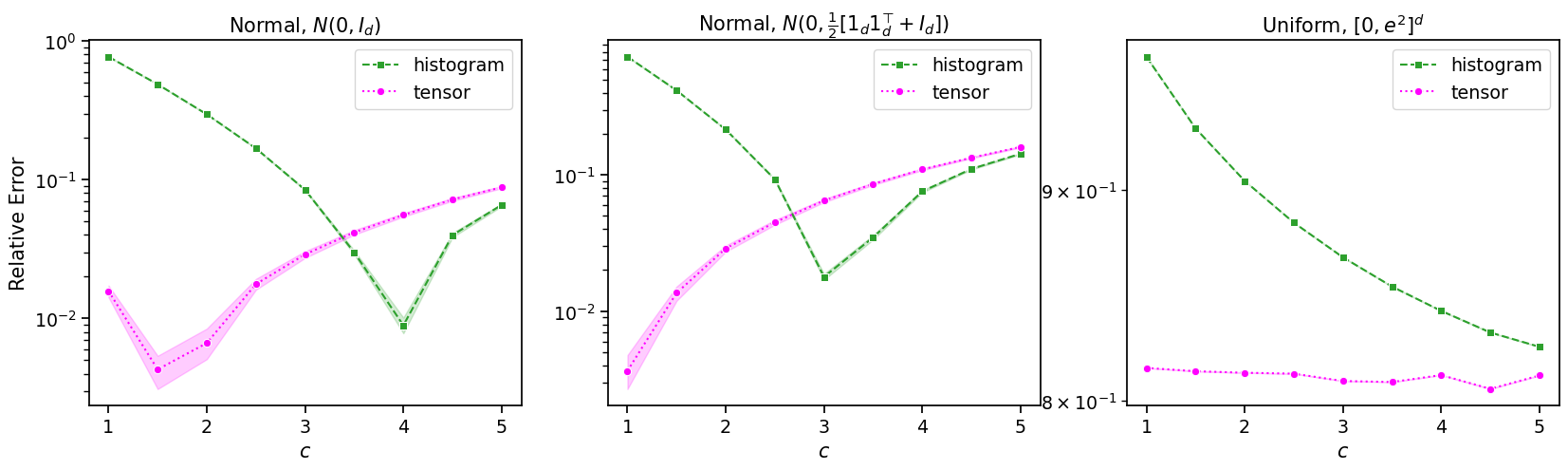}
\caption{Comparing estimates and bins used from 25 trials for dimension $6$ distributions. The histogram is constructed by placing $s=2500$ samples from the distributions into bins of width $c\, s^{-\frac{1}{8}}$ in each dimension, for different values of $c$. Here, the tensor-based approximation uses rank $5$ and the same bins as the histogram-based approximation. 
}\label{fig:bin-schemes}
\end{figure}
While we compare to the $k$-NN estimator for uniform, Gaussian, or $T$ distributions in Figures \ref{fig:tensor-knn-comparisons-uniform} and \ref{fig:tensor-knn-comparisons}, we focus primarily on improvements of PTC over a histogram estimator. 

Preliminary results using PTC and the tensor-based estimator of \cite{vand:21} were comparable in a subset of the experiments of \S\ref{s:bs}, differing in entropy approximation relative error by less than 0.5\%, so we have included only the results of PTC in the experiments presented here. A comprehensive comparison of those two methods is beyond the scope of this work, as PTC can be computed efficiently with sparse data and is guaranteed to be unique (up to scaling and permutation of factors), whereas the other tensor-based approach supports only dense data computations (limiting scalability using sampled histogram tensor data) and is not guaranteed to converge to a unique approximation of the estimator. Further work is needed to address these differences in a fair and complete comparison.

 \subsection{Bin size and samples}\label{s:bs}

In order to determine the influence of bin width, we select two multivariate normal distributions and multivariate uniform distributions with independent variates, each with closed-form formulas for the entropy so that we can determine the relative error in the approximation.
The multivariate normal distribution has entropy $d/2\log(2 \pi e) + 1/2\log \det \Sigma$ where $\Sigma$ is the covariance matrix of order $d$ and the entropy for a multivariate uniform distribution with independent variates is $\log\big((b_1-a_1)\cdots(b_d-a_d)\big)$ where the univariate distribution for the $i$-th variate is over the interval $(a_i,b_i)$. 

\begin{figure}
\centering
\includegraphics[width=0.45\linewidth]{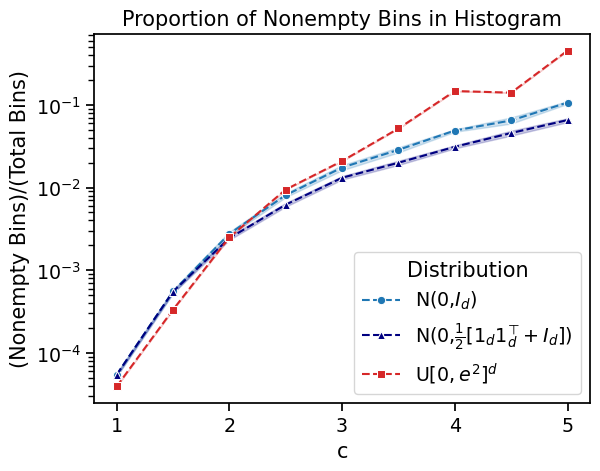}
\caption{The proportion of nonempty bins for the histograms used to estimate entropy in Figure \ref{fig:bin-schemes}. The distributions are dimension $d=6$, and the histograms use bins of width $c\cdot(2500)^{-\frac{1}{8}}$ in each dimension, for different values of $c$. 
The number of bins decreases with increasing $c$. }\label{fig:bin-proportion}
\end{figure}

Figure \ref{fig:bin-schemes} shows that larger bin sizes favor the histogram-based entropy estimates and smaller bin sizes favor the PTC-based entropy estimates, and the difference is nearly two-orders of magnitude for smallest bin sizes.
Figure \ref{fig:bin-proportion} shows the proportion of nonempty bins for the binning schemes in Figure \ref{fig:bin-schemes} with the fraction of nonempty bins occupied differing by four orders of magnitude.
The two figures suggest that increasing the number of bins decreases the fraction of bins with positive counts and so favors the PTC-based entropy estimates, as discussed in \S\ref{subsec:error-analysis}.
The figures also imply that the tensor completion effected has a dramatic impact upon the relative accuracy. 
The three distributions are examples of sub-Gaussian distributions reviewed at the end of \S\ref{sec:tensor-method} and demonstrate that tensor completion favors a binning scheme with small sized bins in contrast to the histogram estimator of the entropy.

\begin{figure}
  \centering 
  \begin{tabular}{c c c}
  \includegraphics[width=.12\linewidth]{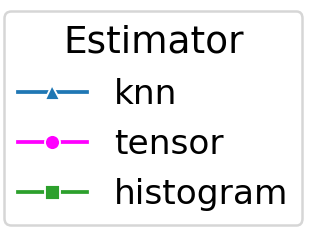}
  &
    \includegraphics[width=0.32\linewidth]{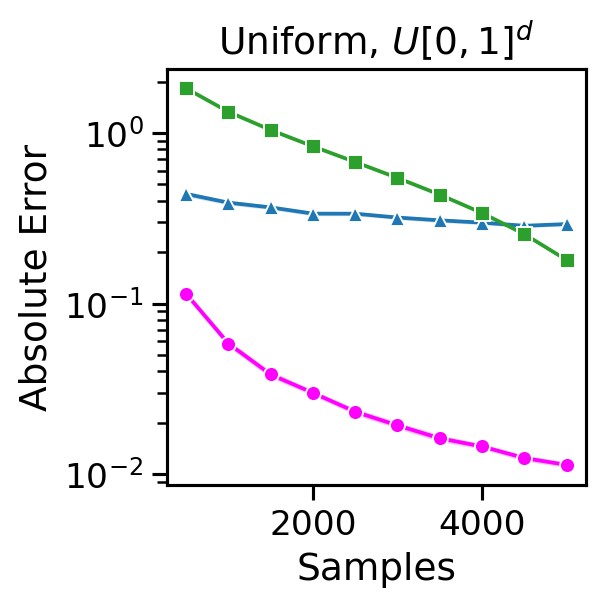}
      &
    \includegraphics[width=0.32\linewidth]{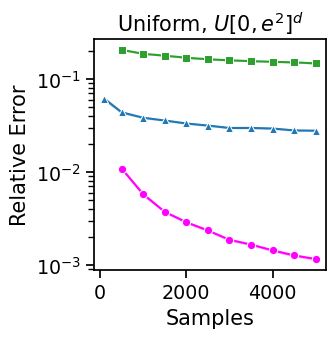}
    \\&\phantom{-------}      (a)&\phantom{------------}       (b)
  \end{tabular}
  \caption{Error in estimated entropy of a five dimensional uniform distribution over (a)  $[0,1]^5$ or (b) $[0,e^2]^5$ with independent dimensions. Estimates use a histogram directly, the tensor (PTC) method, or the $k$-NN method. The results shown are for 25 trials using the $k \in\{1,2,3,\dots,10,25,50,100,200\}$ or rank not exceeding five leading to the smallest error. The histogram uses bins of width $3.5s^{-1/7}$, where $s$ is the number of samples.}
 \label{fig:tensor-knn-comparisons-uniform}
  \end{figure}
  \begin{figure}
  \begin{tabular}{c c c c}
  \includegraphics[width=.10\linewidth]{fixed_dim_tensor_knn_legend.png}
  &
  \includegraphics[width=0.25\linewidth]{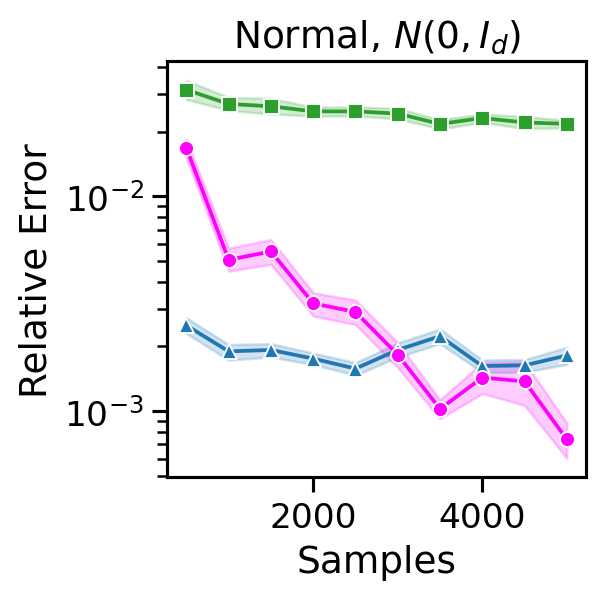}
  &
   \includegraphics[width=0.25\linewidth]{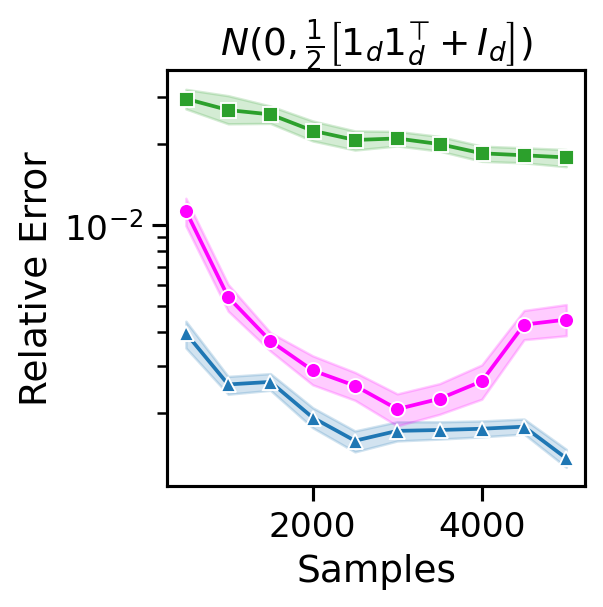}
   &
     \includegraphics[width=0.25\linewidth]{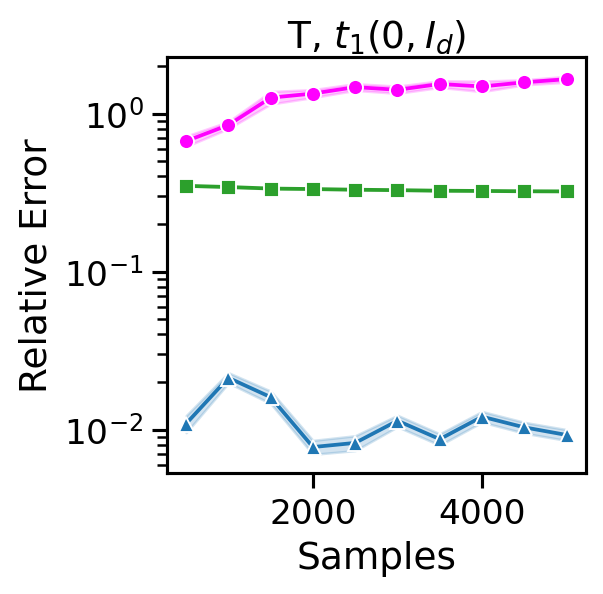}
     \\ &\phantom{-----} (a) &\phantom{-----} (b) &\phantom{-----} (c)
       \end{tabular}
  \caption{Error in estimated entropy using a histogram directly, the tensor (PTC) method, or using the $k$-NN method. The results shown are for dimension 5, over 25 trials using the $k \in\{1,2,3,\dots,10,25,50,100,200\}$ or rank $\leq 5$ leading to the smallest error. The histogram uses bins of width $3.5s^{-1/7}$, where $s$ is the number of samples. The distributions shown are (a) Normal with independent dimensions (b) Normal with correlation between dimensions, and (c) $t$ with one degree of freedom and independent dimensions (equivalent to a Cauchy distribution). }
  \label{fig:tensor-knn-comparisons}
\end{figure}

Figures \ref{fig:tensor-knn-comparisons-uniform} and \ref{fig:tensor-knn-comparisons} show a comparison of the error in estimating differential entropy using a histogram directly, using a low-rank tensor approximation to a histogram, or using the k-NN method for several multivariate distributions with dimension 5. 
In these plots we observe that the PTC estimator outperforms the k-NN estimator for the uniform distribution on $[0,1]^d$ or $[0,e^2]^d$, the two methods are similar (within an order of magnitude) for normal distributions, and the k-NN estimator outperforms the tensor-based estimator for the heavy-tailed Cauchy distribution.
The poor behavior of PTC on the heavy-tailed Cauchy distribution and excellent behavior of PTC on sub-Gaussian distributions suggests that an adequate  number of samples in a small number of bins is important attribute.

\begin{figure}
\centering \includegraphics[width=1.0\linewidth]{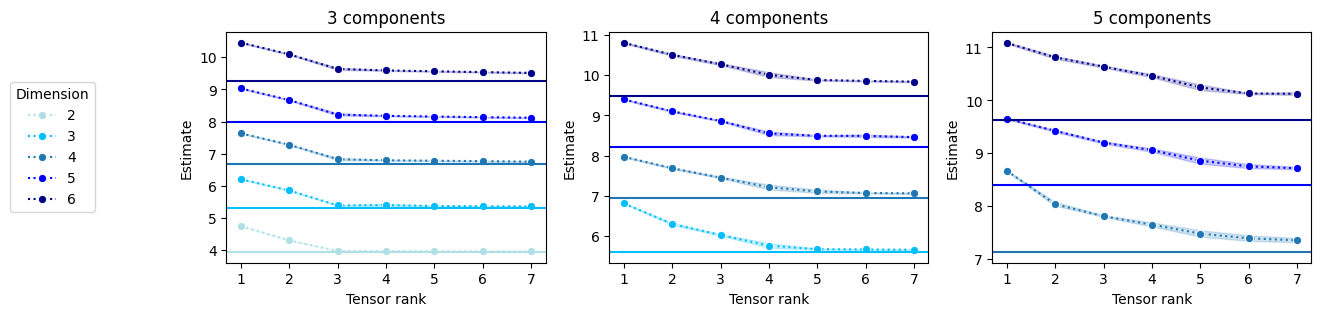}
\caption{Estimated entropy for Gaussian mixtures of dimension $2,3,\dots,6$ with different ranks and different numbers of components with equidistant modes. Dotted lines: tensor estimates from $25$ trials of $s=2500$ samples from the distribution. Solid lines: average histogram estimate in $25$ trials with $s=1000000$ samples. The tensors use 20 bins along each dimension, and the histograms use bins of width $c\cdot s^{-\frac{1}{\text{dim}+2}}$ in each dimension.}
\label{gmm-hist-vs-tens}
\end{figure}

\subsection{Gaussian mixtures}\label{s:gms}
As a step towards considering more general sub-Gaussian distributions, we now investigate the PTC estimator for Gaussian mixtures with different numbers of components. Figure \ref{gmm-hist-vs-tens} shows the estimated entropy for different tensor ranks compared to estimates from histograms with  $s=1,000,000$  samples for three Gaussian mixtures with three, four and five components. For the plots, the Gaussian mixtures were chosen to have components with equidistant modes, each $10$ units apart. Since there are at most $d+1$ equidistant points in $d$-dimensional space, the dimension considered for each number of clusters is at least the number of clusters minus 1.

In the figure, we see a correlation between the number of components in the mixture and the value $R$ at which increasing the tensor rank no longer significantly changes the estimate. The bin sizes for the histogram are given by \eqref{opt-bin-size} and the binning for the tensor estimates uses $20$ bins in each direction. While  there is no closed-form expression for the entropy of a Gaussian mixture with more than one component, the large-sample-size histogram approximation provides a comparison to show that a tensor with a small number of samples can produce an estimate similar to a histogram with a large number of samples. 

\begin{figure}
\centering
\begin{subfigure}[t]{0.22\linewidth}
\centering\includegraphics[width=\linewidth]{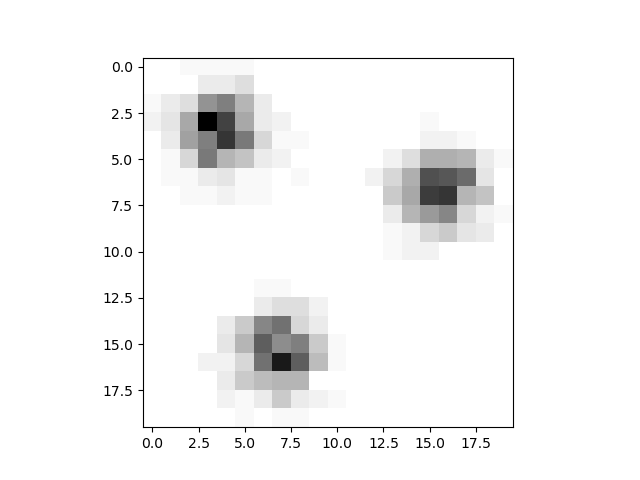}
\caption{Histogram}
\end{subfigure}
\begin{subfigure}[t]{0.22\linewidth}
\centering\includegraphics[width=\linewidth]{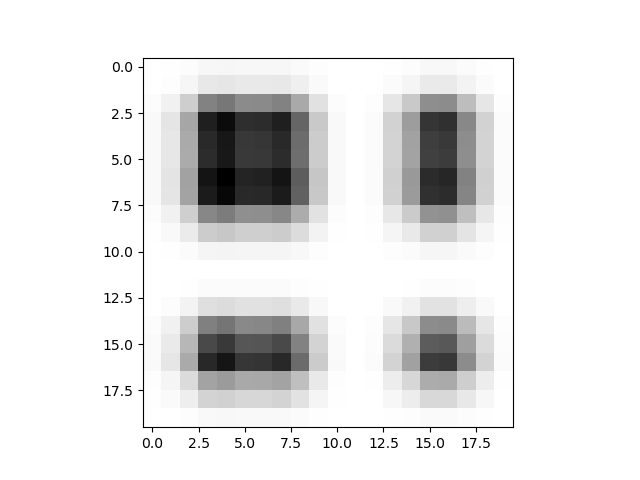}
\caption{Rank 1 decomposition}
\end{subfigure}
\begin{subfigure}[t]{0.22\linewidth}
\centering\includegraphics[width=\linewidth]{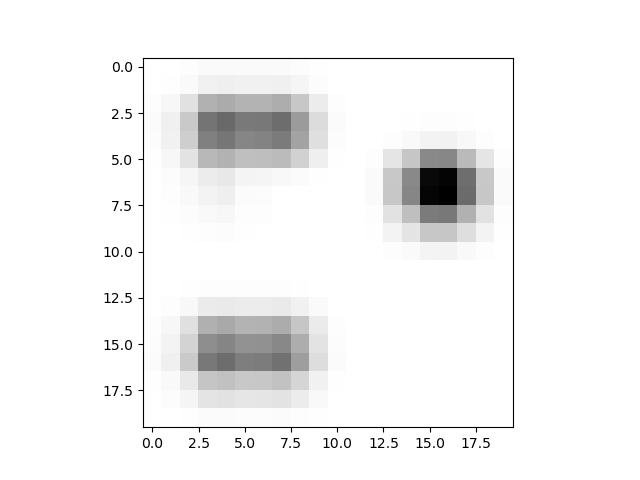}
\caption{Rank 2 decomposition}
\end{subfigure}
\begin{subfigure}[t]{0.22\linewidth}
\centering\includegraphics[width=\linewidth]{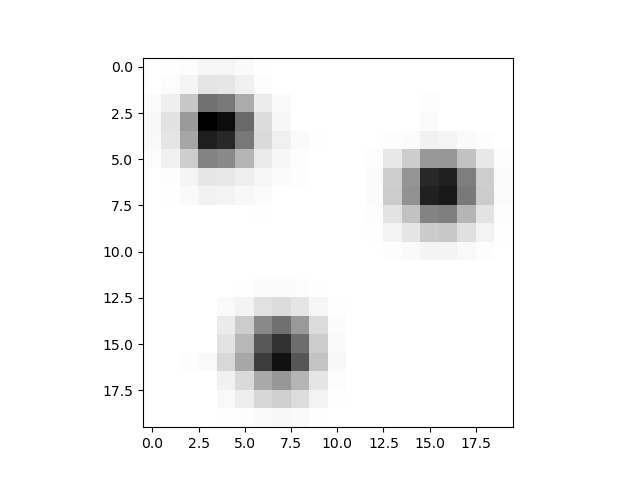}
\caption{Rank 3 decomposition}
\end{subfigure}
\caption{Plots of a histogram of 1000 samples from two-dimensional Gaussian mixture with three equidistant modes and its rank 1, 2, and 3 tensor decompositions.}
\label{fig:2Dtensors}
\end{figure}

A practical detail associated with PTC is the selection of the rank $R$ needed for $\boldsymbol{\widehat{\mathscr{M}}}$.
There is little theory to guide us and determining $R$ depends upon the application, but
the number of components of the mixture provides a clue. 
We observe a relationship between the number of mixture components and the tensor rank used to estimate the entropy. 
Figure \ref{fig:2Dtensors} shows pictures of the tensor decompositions of a histogram for a two-dimensional Gaussian mixture with three components. The sequence suggests that the rank $R$ needs to be as least as large as the number of components.

\begin{figure}
 \centering \includegraphics[width=0.5\linewidth]{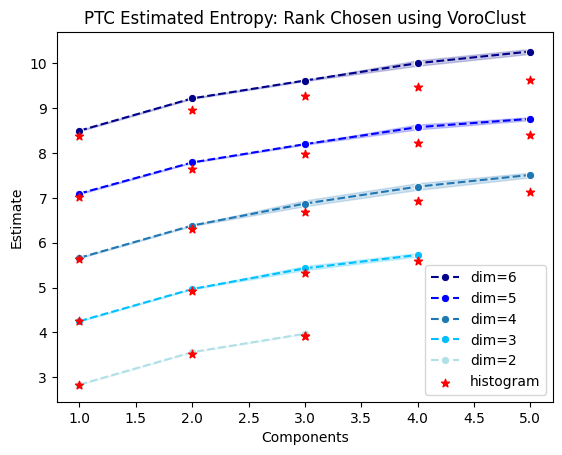}
\caption{Estimated entropy for Gaussian mixtures of dimensions $2,3,\dots,6$ with rank equal to the number of clusters identified by VoroClust and different numbers of components with equidistant modes. The PTC estimates are from $25$ trials, each with $2500$ samples. The $\textcolor{red}{\star}$ points correspond to the average estimate in $25$ trials using a histogram with $s=1000000$ samples and bin widths $3.5\cdot s^{-\frac{1}{\text{dim} + 2}}$ in each dimension. VoroClust was run with adaptively-chosen radii in the sphere cover and allowing $5\%$ of the data to be identified as noise.  The parameters for VoroClust's clustering propagation were  $\texttt{detail\_ceiling}=0.8$ and $\texttt{descent\_limit}=0.1$. In all but 2 trials (both for dimension 2 and 1 component), the number of clusters identified by VoroClust matched the number of components in the Gaussian mixture.}
\label{fig:gmm-with-voroclust}
\end{figure}

Components in mixture models can be viewed as clusters.  This suggests that a clustering tool may be of use for rank selection for a tensor decomposition. There are various clustering tools available. We used VoroClust  \cite{voroclust}, a clustering tool that uses density estimates within a sphere cover of a dataset to identify clusters and noise.  Figure \ref{fig:gmm-with-voroclust} shows the estimated entropy of a Gaussian mixtures using PTC with rank chosen to be the number of clusters identified using VoroClust on the sampled data. The clusters come from a Voronoi tessellation of the dataset's domain. The radii in the sphere cover can be specified or chosen adaptively. We chose to use the adaptive feature to avoid manually tuning radius selection. One benefit of VoroClust is the ability to cluster data with specifying an expected number of clusters. In Figure \ref{gmm-hist-vs-tens}, we observe that the entropy estimate changes slowly with increasing ranks greater than the number of components, so we expect that a small overestimate of the number of clusters would not  cause dramatic difference in the tensor estimates. Although we focused on Gaussian mixtures, we expect that these methods could extend to other mixtures of sub-Gaussian distributions.

 \subsection{Tensor Thresholding} \label{tensor-thres}

The discussion surrounding \eqref{M-def} explains that the factor vectors $\mathbf{a}_r^{(i)}$  can be identified with a probability mass function. 
Hence, the small elements of  $\widehat{\mathbf{a}}_r^{(j)}$ imply that the corresponding elements in the rank-one tensor 
\begin{align*}
\boldsymbol{\widehat{\mathscr{M}}}_r &\coloneqq \widehat{\lambda}_r \, \widehat{\mathbf{a}}_r^{(1)} \circ \widehat{\mathbf{a}}_r^{(2)} \circ \cdots \circ \widehat{\mathbf{a}}_r^{(d)}
\end{align*}
are also small and have negligible contribution to $\text{ent}(\pMh\,\mathbbm{1}_B)$.
For example, if the first element of $ \widehat{\mathbf{a}}_r^{(1)}$ is small, then $n/n_1=n_2 n_3 \cdots n_d$ elements of the $r$-th rank-one tensor of $\boldsymbol{\widehat{\mathscr{M}}}$ are also small. 
This corresponds to a sub-tensor of order $d-1$ for the $r$-th rank-one tensor corresponding to $\widehat{\mathbf{a}}_r^{(1)}$.
This suggests a simple thresholding algorithm that can be used to approximate the entropy estimate in \eqref{ptc-est} and thus lower the memory and computational requirements:
\begin{itemize}
\item
Determine the $R d$ sets $\Omega_{r,i} $ containing the indices for each of the $Rd$ vectors $\widehat{\mathbf{a}}_r^{(i)}$ where the elements  are less than a threshold $0 < \tau < 1$. 
Approximate $\text{ent}(\pMh\,\mathbbm{1}_B)$ using the elements of $\boldsymbol{\widehat{\mathscr{M}}}$ not containing any of the $\sum_{r=1}^R \sum_{i=1}^d \Omega_{r,i}$ indices. 
The number of elements in the sum of the $R$ rank-one tensors is $ R n$. If $|\Omega_{r,i}|$ denotes the number of indices in $\Omega_{r,i}$, then
\begin{align*}
R n - \sum_{r=1}^R \sum_{i=1}^d |\Omega_{r,j}| \frac{n}{n_i} 
\end{align*}
is the number of elements needed to approximate $\text{ent}(\pMh\,\mathbbm{1}_B)$. 
The number of non-negligible elements can be as small as zero when $|\Omega_{r,i}|=n_i$ and as large as $R n$ when $|\Omega_{r,i}|=0$. 
Note that $|\Omega_{r,i}|=n_i$ occurs when $\widehat{\mathbf{a}}_r^{(i)}$ is a uniform distribution (or nearly so) and $\tau n_i > 1$.
\end{itemize}

\begin{figure}
\centering
  \centering \includegraphics[width=.95\linewidth]{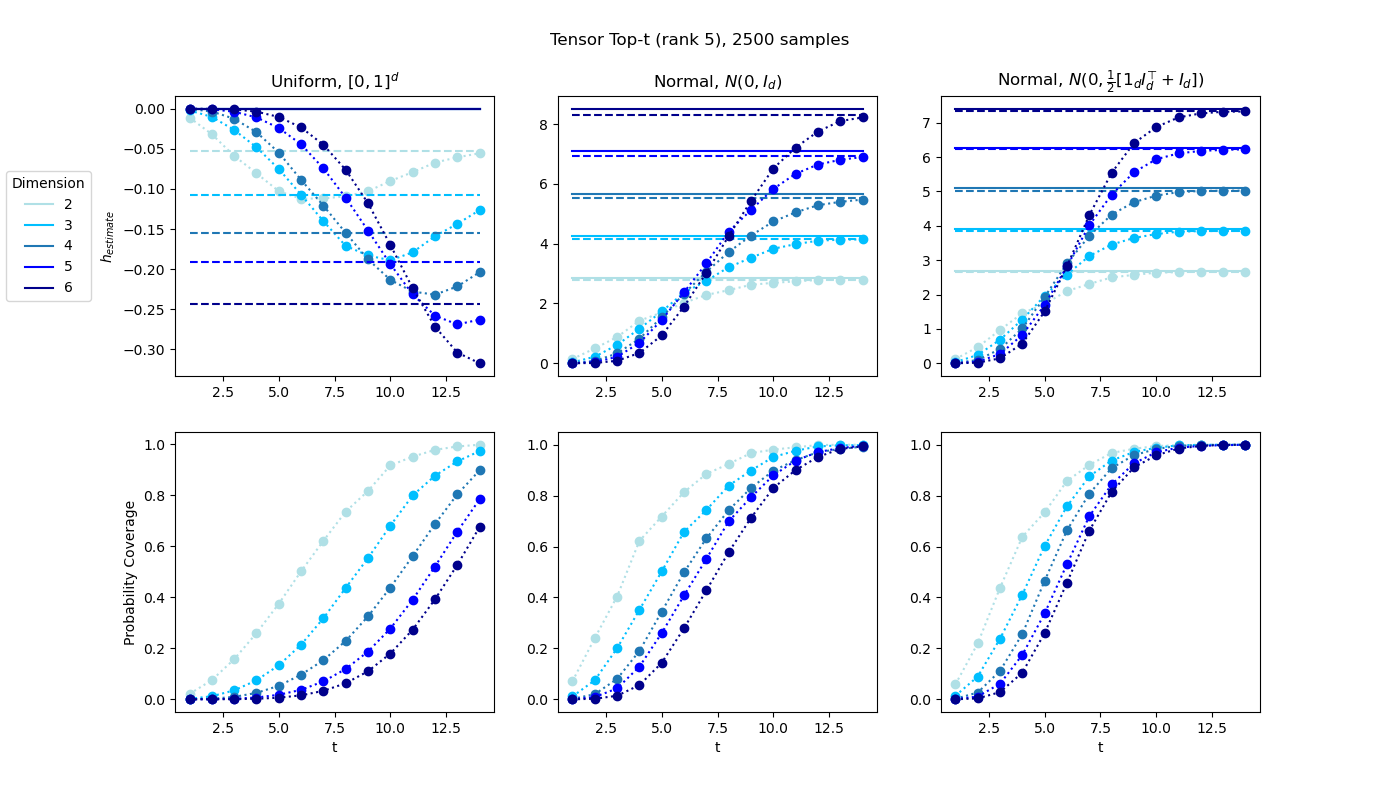}
\caption{The first row shows the entropy estimates using a sample of a tensor, where the indices sampled were the combinations of the $t$ largest indices from each factor vector in (\ref{M-def}) for various $t$. Solid lines show the true entropy and dashed lines show the (full) rank 5 PTC-based estimate. The second row shows the portion of the sum of values in the tensor accounted for through sampling the tensor.}
\label{top-t-plot}
\end{figure}

Figure \ref{top-t-plot} shows entropy estimates found by sampling a PTC estimator  for the uniform distribution on $[0,1]^d$ and normal distributions $N(0,I_d)$ and $N(0,\frac{1}{2}[1_d 1_d^{\top} + I_d])$. 
The estimate to the tensor found by combining the top-$t$ indices approaches the estimate found using the full tensor for increasing $t$.

 \subsection{Application to real-world data}\label{s:real-data}
\begin{figure}
\centering
\includegraphics[width=\linewidth]{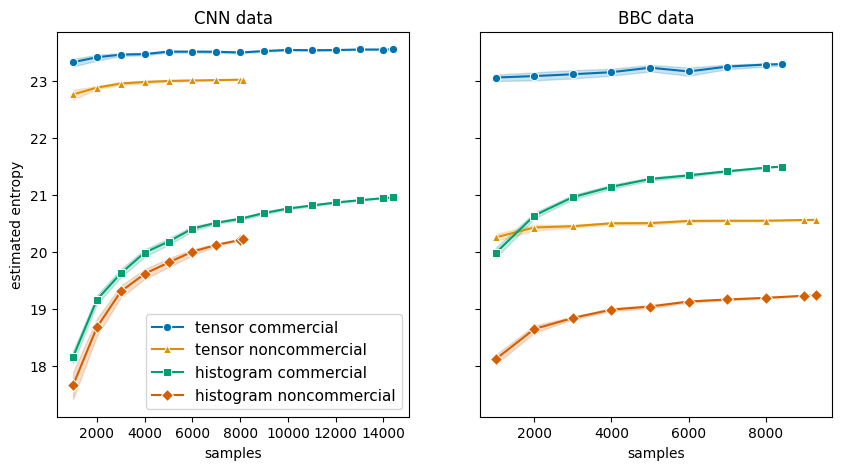}
\caption{Estimated entropy using a PTC estimator or using a histogram. Samples are drawn from the CNN or BBC datasets from \cite{commercial_data}, and we consider 7 features, each identified with a dimension. Both the PTC and histogram estimates use $10$ bins along each dimension. The tensor rank in each trial was chosen using VoroClust with adaptively-chosen radii in the sphere cover and allowing $5\%$ of the data to be identified as noise. The parameters for VoroClust's clustering propagation were  $\texttt{detail\_ceiling}=0.8$ and $\texttt{descent\_limit}=0.1$. Samples were drawn without replacement, and the points in the plot show the average estimate in 25 trials. The last point for each line uses every available sample from the corresponding dataset.}
\label{fig:TVdata}
\end{figure}

To demonstrate the PTC estimator on real-world data, we approximate the differential entropy using 7 features (mean values of frame difference distribution, short time energy, ZCR, spectral centroid, spectral roll off, spectral flux, and fundamental frequency) in the CNN and BBC broadcast news datasets from \cite{commercial_data}.  We identify each of the selected features with a variate. To reduce the memory requirements for the tensors and histograms, we used $10$ bins along each dimension ($10^7$ total) rather than $20$ bins along each dimension as in our previous experiments. While we cannot compare to the (unknown) true entropy value for the distribution, we can compare PTC and histogram-based estimates of entropy as the number of samples increases, as shown in Figure \ref{fig:TVdata}. The datapoints are each labeled as ``commercial'' or ``noncommercial.''  We observe that the changes in the PTC-based estimates are smaller than the changes in the histogram-based estimates as the number of samples increases, and in the CNN data, the PTC-based estimate distinguishes between sets of ``commercial'' and ``noncommercial'' points with smaller sample sizes. Using all available points, the CNN data histograms were significantly more sparse than the CNN data tensors (approximately  $ 99.96\%$ zeros compared to an average of approximately $20\%$ zeros). The BBC data histograms using all available data were approximately $ 99.98\%$ zeros, and the BBC data tensors were on average approximately $83.2\%$ zeros. In most trials, VoroClust identified one cluster in the CNN data and two clusters in the BBC data.

\section{Conclusion}
\label{sec:conclusion}

We introduced the \emph{Poisson tensor completion} (PTC) parametric estimator 
to approximate the density of a multivariate distribution.
The PTC estimator exploits the inter-sample relationships determined during a maximum Poisson likelihood estimation introduced in~\cite{chikolda12} by Chi and Kolda to compute a low-rank Poisson CP tensor decomposition. 
Our crucial observation is that the histogram bins are an instance of a space partitioning of counts and thus can be identified with a spatial non-homogeneous Poisson process, which leads to the density estimate.
Such an identification and application to density estimation is novel with this paper.
The Poisson tensor decomposition leads to a completion of the mean measure over all bins---including those containing few to no samples---and leads to our proposed PTC estimator.

Our error analysis underscored the role played by the tensor completion to impute values for the density and demonstrate the impact of our estimator approach on computing expectations such as differential entropy. 
Several numerical experiments on synthetic and real-world data examined the role of sub-Gaussian distributions, thresholding, and tensor rank on the quality of the PTC estimator compared with a histogram estimator as measured by using each as a plug-in estimator when computing differential entropy.
PTC works well on sub-Gaussian distributions when there is an adequate number of samples in a small number of bins.
We also explained that heavy-tailed distributions, which do not concentrate in norm, do not benefit from the PTC estimator. 
Our experiments also indicate that in practice, the number of variates is limited to available computational resources and  further improvements achievable by adroit thresholding.

We reviewed recent work describing the role of the tensor rank in  the identifiability and indeterminacy of the PTC parametric estimator. 
The error analysis also quantified the effect of tensor rank and the numerical experiments investigated the role played by tensor rank. 
We demonstrated that the number of components in a mixture model is correlated with the tensor rank and can be determined via the use of clustering software.

Our future work is a rigorous  analysis for the approximation provided by tensor completion and appropriate binning strategies. 
Previous work on Poisson tensor completion demonstrated that a zero-truncated Poisson CP decomposition~\cite{lopez23} can better estimate expected counts than when using a Poisson CP decomposition if the number of zero counts in the observed tensor is not too large relative to the sizes of the tensor dimensions. 
We plan to explore the conditions where such a result combined with the thresholding approach above may lead to an even more efficient estimator for sub-Gaussian distributions based on zero-truncated Poisson tensor completion.  

\section{Acknowledgments}
We thank Carlos Llosa and Derek Tucker  of Sandia National Labs, and Oscar Lopez of Florida Atlantic University for several helpful discussions during the writing of this manuscript,  
The second author thanks Scott McKinley of Tulane University for discussions on spatial Poisson processes.
Three referees provided constructive reviews that lead to a broadening of the original scope of our paper. 
In particular, one of the referees made the astute observation that our  estimator is parametric instead of the originally proposed non-parametric estimator.

This work was supported by the Laboratory Directed Research and Development program (Project 233076) at Sandia National Laboratories, a multimission laboratory managed and
operated by National Technology and Engineering Solutions of Sandia LLC, a wholly owned subsidiary of Honeywell International Inc. for the U.S. Department of Energy’s National Nuclear Security Administration under contract DE-NA0003525.
\bibliographystyle{siamplain}
\bibliography{multivariate-entropy}

\end{document}